\documentclass{article}
\usepackage{epsfig}
\usepackage{psfrag}

\usepackage{amsmath,amssymb,amsthm}

\usepackage{hyperref}

\title{
Curious Continued Fractions, Nonlinear Recurrences and Transcendental Numbers}
\author{Andrew Hone\thanks{School of Mathematics,
Statistics and Actuarial Science, University of
Kent, Canterbury CT2 7NF, U.K. ~~E-mail: A.N.W.Hone@kent.ac.uk}
}

\newcommand{\beq}{\begin{equation}}
\newcommand{\eeq}{\end{equation}}
\newcommand{\bear}{\begin{array}}
\newcommand{\eear}{\end{array}}
\newcommand\la{{\lambda}}
\newcommand\La{{\Lambda}}
\newcommand\ka{{\kappa}}
\newcommand\al{{\alpha}}

\newtheorem{thm}{Theorem}[section]

\newtheorem{propn}[thm]{Proposition}

\newtheorem{rem}[thm]{Remark}

\newtheorem{cor}[thm]{Corollary}

\theoremstyle{remark}

\newcommand{\Z}{{\mathbb Z}}

\newcommand{\seqnum}[1]{\href{http://oeis.org/#1}{\underline{#1}}}

\begin{document}

\maketitle

\begin{abstract}
We consider a family of integer sequences generated by nonlinear recurrences of the second order, which have 
the curious property that the terms of the sequence, 
and integer multiples of the ratios of successive terms (which are also integers), 
appear interlaced in the continued fraction expansion of the sum of the reciprocals of the terms. Using the rapid 
(double exponential) growth 
of the terms, for each sequence it is shown that the sum of the reciprocals is a transcendental number. 
\end{abstract}

\section{Introduction}

\setcounter{equation}{0}

For some time there has been considerable interest in rational recurrences which surprisingly generate integer sequences. 
Such sequences were made popular by the articles of Gale \cite{gale,gale2}, who discussed 
some particular nonlinear recurrence 
relations of the form 
\beq 
\label{recf} 
x_{n+N} \, x_n = f(x_{n+1}, \ldots ,  x_{n+N-1}), 
\eeq 
where $f$ is a polynomial in $N-1$ variables. Observe that the above recurrence is rational, in the sense 
that each new iterate $x_{n+N}$ is a rational function of the $N$ previous terms $x_n,\ldots,x_{n+N-1}$. Starting from 
$N$ initial values $x_0,\ldots,x_{N-1}$ which are all integers, there is no reason to expect that subsequent terms will be, because 
one must divide by $x_n$ at each step. However, a very wide variety of examples are now known, for which the recurrence (\ref{recf}) has the Laurent 
property: if the initial values are viewed as variables, then for certain special  choices of $f$, all of the iterates belong to the ring $\Z[x_0^{\pm 1},\ldots ,x_{N-1}^{\pm 1}]$, consisting of 
Laurent polynomials in the initial values with integer coefficients. In particular, the Laurent property implies that if all the initial values are taken 
to be 1 (or $\pm 1$), then $x_n\in\Z$ for all $n$. 

The Laurent property is a key feature of Fomin and Zelevinsky's cluster algebras \cite{fz1}, which are generated by birational iterations of the same 
shape as (\ref{recf}), that is 
$$ 
\mathrm{old}\,\,\mathrm{variable}\times \mathrm{new}\,\,\mathrm{variable}=\mathrm{exchange}\,\,\mathrm{polynomial},
$$
in the particular case that the exchange polynomial $f$ is a binomial.  
The main tool available for proving the Laurent property is the Caterpillar Lemma 
due to Fomin and Zelevinsky \cite{fz}, which also applies to more general choices of $f$, 
fitting into the broader framework of Laurent Phenomenon (LP) algebras \cite{lp}. Within the axiomatic setting of cluster algebras or LP algebras, 
there is a 
requirement that the exchange polynomials should not be divisible by any of the variables. However, this requirement is not necessary for the Laurent 
property to hold. Indeed, even for the case of a recurrence of second order, of the form 
\beq\label{rec2}
x_{n+2} \, x_n = f(x_{n+1}), 
\eeq 
the requirement that $x\not |\, f(x)$ is not necessary. In work by the author \cite{honepla,numberpoly} it was shown that recurrences of the form 
(\ref{rec2}) having the Laurent property fit into three classes, depending on the form of $f$: (i) $f(0)\neq 0$, in which case the recurrence belongs 
within the framework of cluster algebras (when it is a binomial) or LP algebras (when it is not); 
(ii)  $f(0)= 0$, $f'(0)\neq 0$; (iii) $f(0)=f'(0)= 0$. 
In classes (i) and (ii) there are additional requirements on $f$, but in class (iii) one can take $f(x)=x^2 F(x)$ with arbitrary $F\in\Z[x]$. 

The simplest non-trivial example of the form  (\ref{rec2}) belonging to the third of the classes identified in 
previous work by the author \cite{honepla} is the recurrence 
\beq\label{orig} 
 x_{n+2} \, x_n = x_{n+1}^2(x_{n+1}+1). 
\eeq 
Due to the Laurent property, the initial values $x_0=x_1=1$ generate an integer sequence: 
\beq\label{exs} 
1,1,2,12,936,68408496,342022190843338960032,\ldots;
\eeq 
this is sequence \seqnum{A112373} in the Online Encyclopedia of Integer Sequences (OEIS). 
As one might expect from the first few terms, this sequence grows very rapidly: $\log x_n \sim C\la^n$ with $C\approx 0.146864$ and 
$\la=(3+\sqrt{5})/2$. Another feature of sequence \seqnum{A112373} is that the 
ratios $y_n=x_{n+1}/x_n$ also form an integer sequence, that is 
\beq\label{ratios}
1,2,6, 78,73086,4999703411742,1710009514450915230711940280907486,\ldots, 
\eeq 
which is sequence  \seqnum{A114552} in the OEIS, and the same is true for the ratios of ratios, i.e., $z_n=y_{n+1}/y_n=x_{n+1}+1$ by (\ref{orig}); 
this property of the ratios is common to all recurrences in class (iii).  
Hanna made some very interesting empirical observations about the sequence (\ref{exs}) \cite{hanna}, by considering $\mathcal{S}$, the sum of reciprocals of the terms:  
\beq\label{sumr} 
\mathcal{S}=\sum_{j=0}^\infty \frac{1}{x_j} = 1+1+\frac{1}{2}+\frac{1}{12}+\frac{1}{936}+\cdots \approx 2.5844017240.
\eeq  
In the OEIS, the digits of this number appear as  sequence  \seqnum{A114550}, yet it is not the decimal expansion of $\mathcal{S}$ that is interesting, but 
rather its continued fraction representation; with the notation 
$$ 
[a_0; a_1,a_2,a_3,\ldots,a_n,\ldots ]=a_0+\cfrac{1}{a_1+\cfrac{1}{a_2+\cfrac{1}{a_3+\cdots\cfrac{1}{a_n+\cdots}}}}, 
$$ 
one finds that 
\beq\label{cfrs} 
\mathcal{S}=[2;1,1,2,2,6,12,78,936,73086,68408496,4999703411742,\ldots].
\eeq 
What appears to be the case from the above is that (apart from the initial value $a_0=2$) 
the sequence of $a_n$ in (\ref{cfrs}), which is 
number  \seqnum{A114551}, is obtained by  
interlacing the original sequence (\ref{exs}) with the ratios (\ref{ratios}). 
As observed by Shallit \cite{shallit}, and shown by Harris \cite{harris}, 
this implies that the even/odd 
terms satisfy 
\beq\label{shallit} 
a_{2n}=a_{2n-1}a_{2n-2}, \qquad  a_{2n+1}=a_{2n-1}(a_{2n}+1), 
\eeq 
respectively, where the first formula holds for $n\geq 2$ and the second for $n\geq1$.  

The purpose of this short note is to prove Hanna's observations concerning the continued fraction expansion (\ref{cfrs}), and generalize 
them to an infinite family of integer sequences generated by recurrences belonging to class (iii) in the author's 
previous work \cite{honepla}. At the same time we show that the 
number  $\mathcal{S}$ given by (\ref{sumr}) is transcendental, and the same is true for the sums of reciprocals obtained from the other sequences in this 
family. 

The results presented here  are similar in spirit to those in a  paper by Davison and Shallit, who found some continued fractions 
whose partial quotients (coefficients) are explicitly related to the denominators of  their convergents, and used this to prove the 
transcendence of Cahen's constant \cite{ds}. For references to examples of other transcendental numbers whose 
complete continued fraction expansion is known, the reader should consult the latter paper. 

\section{Continued fractions for sums of reciprocals} 

\setcounter{equation}{0}

Before we proceed with presenting a family of sequences which generalizes  (\ref{exs}), we present some  facts about this particular example, to 
motivate the proof of the main result. When taking the sum of reciprocals  (\ref{sumr}), it is more convenient to 
exclude the index $j=0$ from the sum, and then consider the partial sums 
$$
S_N=\sum_{j=1}^N \frac{1}{x_j}. 
$$ 
Calculating the finite continued fractions of these partial sums, we find that
$$ 
S_1=1,\quad S_2=\frac{3}{2}=1+\frac{1}{2},\quad S_3=1+ \cfrac{1}{1+\cfrac{1}{1+\cfrac{1}{2+\cfrac{1}{2}}}},
$$ 
$ 
S_4=[1;1,1,2,2,6,12]$, 
$S_5=[1;1,1,2,2,6,12,78,936]$, and 
$$
S_6=[1;1,1,2,2,6,12,78,936,73086,68408496]
$$ 
are the first few partial sums. As will be proved in due course, the pattern is 
\beq\label{pat}
S_N=[x_0; y_0,x_1,y_1,x_2,\ldots, y_{N-2},x_{N-1}], 
\eeq
so that the even/odd coefficients are $a_{2n}=x_n$ and $a_{2n+1}=y_n$ respectively, and as we have chosen to start the sum with $1/x_1=1$ 
we now have $a_0=x_0=1$ which ensures that both formulae (\ref{shallit}) hold for all $n\geq 1$. The result for $S_2$ looks anomalous, but in fact 
(\ref{pat}) is seen to hold for $N=2$ upon noting that 
$$ 
S_2=[1;2]=[1;1,1].
$$ 
The continued fraction for the infinite 
sum $S_\infty=\sum_{j=1}^\infty \frac{1}{x_j}$ is obtained in the limit $N\to\infty$, and compared with (\ref{sumr}) we have 
$\mathcal{S}=S_\infty+1$. 

We now wish to generalize these observations to integer sequences generated by recurrences of the 
form 
\beq\label{form} 
x_{n+2} \, x_n = x_{n+1}^2\, F(x_{n+1}), 
\eeq 
where $F(x)\in\Z[x]$, and we assume that $d=\deg F\geq 1$ to avoid a trivial case. 
If we take such a recurrence with the initial values $x_0=x_1=1$, and set 
$$
y_n=\frac{x_{n+1}}{x_n}, \qquad z_n=\frac{y_{n+1}}{y_n}=\frac{x_{n+2}x_n}{x_{n+1}^2}=F(x_{n+1}), 
$$ 
then we have $y_0=1$, $x_2=y_1=z_0=F(1)$, and by induction we see that $x_n,y_n,z_n\in\Z$ for all $n\geq 0$, 
so we have three integer sequences, 
as long as 
the recurrence (\ref{form}) does not reach a singularity (division by zero), 
which can happen if $F(x)=0$ for some $x\in\Z$. Indeed, suppose that for some positive integer $m$ we have 
$0\neq x_n\in\Z$ for $0\leq n\leq m$, but $F(x_m)=0$; then at the next step $x_{m+1}=0$, followed by $x_{m+2}=0$, 
and then $x_{m+3}$ is undefined. The analysis of these recurrences near to a singularity has been performed  
previously \cite{honepla}.
 
To avoid the possibility of reaching a singularity from the initial conditions $x_0=x_1=1$,
we choose $F$ to have only positive integer coefficients ($F(x)\in\Z_{\geq 0}[x]$), 
so that $F(x)>0$ whenever $x>0$, and then 
all three sequences, $(x_n)$, $(y_n)$ and $(z_n)$, consist of positive integers. In order 
for the continued fraction expansion of the sum of reciprocals to behave in the right way, 
we must make the further assumption that $F(0)=1$, 
which (since the degree of $F$ is positive) implies $F(1)>1$, 
and hence $x_{n+1}>x_n$ for $n\geq 1$, and the sequences $(y_n)$ and 
$(z_n)$ are strictly increasing as well. The precise rate of growth of these sequences will be 
considered in the next section, but for now 
we proceed with the main result on continued fractions. 
\begin{thm}\label{maint}
For a sequence $(x_n)$ generated from the initial values $x_0=x_1=1$ 
by the recurrence (\ref{form}) with $F(x)\in\Z_{\geq 0}[x]$ and $F(0)=1$, the partial sums of reciprocals have the continued fraction expansions 
\beq\label{parcf} 
S_N=\sum_{j=1}^N \frac{1}{x_j} =[a_0;a_1,a_2,\ldots,a_{2N-2}]
\eeq 
for all $N\geq 1$, where 
\beq \label{coeffs} 
a_{2n}=x_n, \qquad a_{2n+1}=\frac{F(x_{n+1})-1}{x_n}\in\Z_{>0}. 
\eeq
\end{thm} 
\begin{cor}\label{infs} 
The infinite sum of the reciprocals of the terms of the sequence $(x_n)$ is given by 
\beq\label{infcf} 
S_\infty=\sum_{j=1}^\infty \frac{1}{x_j} =[a_0;a_1,a_2,\ldots,a_n, \ldots],
\eeq 
where the coefficients 
of the continued fraction are as in (\ref{coeffs}).
\end{cor} 
The latter result on the infinite sum follows immediately from (\ref{parcf}) by taking the limit $N\to\infty$. Note that the  
convergence of the sum is guaranteed since the infinite continued fraction makes sense for any sequence of positive integer 
coefficients $(a_n)$; the convergence of the infinite sum can also be proved directly by using estimates for the growth 
of $(x_n)$, as given in the next section. 
\begin{rem} 
We may write 
$$
S_\infty-1=\sum_{j=1}^\infty \frac{1}{y_1y_2\cdots y_j}
$$ 
and observe that $(y_n)$ is a non-decreasing sequence of positive integers, with $y_n\geq 2$ for $n\geq 1$, so this is 
an example of an Engel expansion (see Theorem 2.3 in Duverney's book \cite{duve}). 
\end{rem} 

To prove Theorem \ref{maint}, we start by recalling some general facts about convergents of continued fractions, which are well known; for a 
brief summary of these results, the reader is referred to the first chapter of Manin and Panchishkin's 
book \cite{manin}, or for more details see the book by Cassels \cite{cassels}. The 
$n$th convergent of the continued fraction $[a_0;a_1,a_2,\ldots]$ is given by 
$$ 
\frac{p_n}{q_n} = [a_0;a_1,a_2,\ldots,a_n], 
$$ 
where the numerators $p_n$ and denominators $q_n$ are given in terms of the coefficients according to the matrix identity 
\beq\label{matid} 
\left(\begin{array}{cc} a_0 & 1 \\ 1 & 0 \end{array}\right)
\left(\begin{array}{cc} a_1 & 1 \\ 1 & 0 \end{array}\right) 
\ldots 
\left(\begin{array}{cc} a_n & 1 \\ 1 & 0 \end{array}\right) = 
\left(\begin{array}{cc} p_n & p_{n-1} \\ q_n & q_{n-1} \end{array}\right), 
\eeq 
and both $p_n$ and $q_n$ are obtained recursively via the same linear three-term recurrence relation, that 
is 
\beq\label{3term} 
\begin{array}{rcl} 
p_{n+1}& = & a_{n+1}p_n+p_{n-1}, \\ 
q_{n+1}& = & a_{n+1}q_n+q_{n-1},
\end{array} 
\eeq 
with the initial values 
\beq\label{inits} 
q_{-1}=0, \qquad p_{-1}=q_0=1, \qquad p_0=a_0.   
\eeq   
Taking the determinant of both sides of (\ref{matid}) gives the formula 
\beq\label{detid} 
p_nq_{n-1}-p_{n-1}q_n = (-1)^{n+1}, 
\eeq 
valid for $n\geq 0$,
and the identity 
\beq\label{nid} 
p_nq_{n-2}-p_{n-2}q_n =  (-1)^{n}a_n
\eeq 
for $n\geq 1$ follows by combining (\ref{3term}) with (\ref{detid}). 

Now consider the 
convergents of the continued fraction whose coefficients $a_n$ are given in terms of the sequence 
$(x_n)$ by (\ref{coeffs}). First of all, note that if we write $F(x)=1+xG(x)$ for $G\in\Z_{\geq 0}[x]$, then 
$a_{2n+1}=x_{n+1}G(x_{n+1})/x_n=y_nG(x_{n+1})\in\Z_{>0}$ as claimed; so in general the coefficients with odd index 
are integer multiples of the ratios $y_n$. We prove by induction that  
\beq\label{indn} 
q_{2N-1}=y_N-1=\frac{x_{N+1}}{x_N}-1, \qquad q_{2N} =x_{N+1}. 
\eeq 
For $N=0$ we have $q_{-1}=x_1/x_0-1=0$ and $q_0=x_1=1$ in agreement with (\ref{inits}), and assuming 
that (\ref{indn}) holds for some $N$, from (\ref{3term}) we have 
$$ 
q_{2N+1}=a_{2N+1}q_{2N}+q_{2N-1}=(F(x_{N+1})-1)\frac{x_{N+1}}{x_N}+\frac{x_{N+1}}{x_N}-1
=\frac{x_{N+2}}{x_{N+1}}-1, 
$$  
by (\ref{form}), and hence  
$$
q_{2N+2} =  a_{2N+2}q_{2N+1}+q_{2N}=x_{N+1}\left(\frac{x_{N+2}}{x_{N+1}}-1\right)+x_{N+1} =
x_{N+2}
$$
as required. Now  (\ref{parcf}) is clearly true for $N=1$, and if it holds for some index $N$ then, by using the second 
equation in (\ref{indn}) as well as (\ref{nid}) and (\ref{coeffs}), 
we find 
$$ 
S_{N+1}= S_N+\frac{1}{x_{N+1}}=\frac{p_{2N-2}}{q_{2N-2}}+\frac{1}{q_{2N}} =\frac{p_{2N}q_{2N-2}-a_{2N}+q_{2N-2}}{ q_{2N-2}q_{2N}} 
=\frac{p_{2N}}{q_{2N}}
$$ 
which means that the partial sum  $S_{N+1}$ is the $2N$th convergent of the continued fraction $[a_0;a_1,a_2,\ldots]$ with 
coefficients (\ref{coeffs}), so (\ref{parcf}) holds for the index $N+1$. 
This completes the proof of Theorem \ref{maint}. 

\section{Transcendence of the sums} 

\setcounter{equation}{0} 

We now prove the following 
\begin{thm}\label{trans} 
The infinite sum $S_\infty$, as in  (\ref{infcf}), is a transcendental number.   
\end{thm} 
The proof is based on Roth's theorem, which says that if $\al$ is an irrational algebraic number then for an arbitrary fixed $\delta>0$ there are only 
finitely many rational approximations $p/q$ for which 
\beq\label{roth}
\left|\al -\frac{p}{q}\right|<\frac{1}{q^{2+\delta}} 
\eeq 
(see Manin and Panchishkin's 
book for a brief discussion \cite{manin}, and for a proof see the book by Cassels \cite{cassels}).  
Note that the number $S_\infty$ is irrational, since its continued fraction expansion (\ref{infcf}) 
consists of an infinite sequence of coefficients $a_n\neq0$. 
In order to make use of Roth's theorem, it is enough for the integer sequence 
$(x_n)$ to satisfy the growth condition 
\beq\label{growth} 
x_{n+1}>x_n^\ka \qquad \mathrm{for}\,\,\mathrm{some}\,\,\ka>2, 
\eeq 
for all sufficiently large $n$. Supposing that this is so, it follows that 
$$ 
x_{n+j}>x_n^{\ka^j} \qquad \mathrm{for} \,\,j\geq 1 
$$ 
whenever $n$ is large enough. 
Then from Theorem \ref{maint} we have 
$$ 
\left|S_\infty-\frac{p_{2n}}{q_{2n}}\right|=\sum_{j=n+2}^\infty \frac{1}{x_j} <\sum_{j=1}^\infty \frac{1}{x_{n+1}^{\ka^j}}. 
$$   
Now the function $g(j)=j^{\frac{1}{j-1}}$ is monotone decreasing with $g(2)=2$, so $\ka^j>j\ka$ for $j\geq 2$ and $\ka>2$, which 
together with (\ref{indn}) implies 
$$
\left|S_\infty-\frac{p_{2n}}{q_{2n}}\right|< \sum_{j=1}^\infty \frac{1}{x_{n+1}^{j\ka}}=\frac{(1-x_{n+1}^{-\ka})^{-1}}{x_{n+1}^{\ka}} 
< \frac{1}{x_{n+1}^{\ka-\epsilon}}=\frac{1}{q_{2n}^{\ka-\epsilon}}
$$ 
for any $\epsilon>0$ and $n$ sufficiently large. So if $\epsilon$ is chosen such that $\ka-\epsilon=2+\delta >2$, then $\al=S_\infty$ has infinitely many 
rational  approximations satisfying (\ref{roth}), and hence must be transcendental. 

To show that (\ref{growth}) holds for any sequence $(x_n)$ defined by the recurrence (\ref{form}) with $x_0=x_1=1$ and $F(0)=1$, $F(x)\in\Z_{\geq 0}[x]$
as in Theorem \ref{maint}, we can use a very crude estimate. Indeed, we have $F(x)=1+\ldots +cx^d$ for some integer $c\geq1$, 
so $F(x)>x^d$ for all $x>0$, and then from  (\ref{form}) we obtain 
\beq \label{vcrude} 
x_{n+1}>\frac{x_n^{d+2}}{x_{n-1}}  \geq x_n^{d+1} 
\eeq 
for all $n\geq 1$, since the ratios $y_n=x_{n+1}/x_n$  form an increasing sequence. The above 
growth condition is sufficient for (\ref{growth}) when the degree $d\geq2$, but not when $d=1$, which is 
the case relevant to the original sequence (\ref{exs}). However, this estimate can be improved upon 
by using $x_{n-1}<x_{n}^{\frac{1}{d+1}}$ in the first inequality in (\ref{vcrude}), to yield
\beq \label{crude} 
x_{n+1}>x_n^{d+2-\frac{1}{d+1}}\geq x_n^{5/2}  
\eeq 
for $d\geq 1$. 

In fact, we can get a more accurate measure of growth from asymptotic arguments. Upon 
taking the logarithm of (\ref{form}) we find that $\La_n=\log x_n$ satisfies 
\beq\label{lineq} 
\La_{n+1}-(d+2)\La_{n}+\La_{n-1} =\log c +\al_n, \quad \mathrm{with}\quad 
\al_n =\log\left(\frac{F(x_n)}{cx_n^d}\right). 
\eeq 
Note that $\al_n =\log(1+O(x_n^{-1}))=O(x_n^{-1})$ as $n\to\infty$, 
which means that, to leading order, the growth of $\La_n$ is governed by a homogeneous 
linear equation with constant coefficients, given by the vanishing of the left-hand side of (\ref{lineq}). 
The characteristic  equation is $\la^2-(d+2)\la +1=0$, with the largest root being  
\beq\label{root} 
\la=\frac{d+2+\sqrt{d(d+4)}}{2}>2 .
\eeq 
Thus we find that   
\beq\label{asymp} 
\La_n\sim C\la^n,  
\eeq 
for some $C>0$. Hence $\La_{n+1}/\La_n\to \la$ as $n\to\infty$, and so for any $\epsilon>0$ it follows that 
$\La_{n+1}>(\la-\epsilon)\La_n$ for all sufficiently large $n$, giving the required growth condition 
$$ 
x_{n+1}>x_n^{ \la-\epsilon}. 
$$  

The asymptotic properties of the sequence $(x_n)$ can be determined more precisely by adapting the methods 
of Aho and Sloane \cite{ahos}, leading to the following result, which is easy to verify directly from (\ref{lineq}). 
\begin{propn} For the initial conditions $x_0=x_1=1$, the logarithm $\Lambda_n=\log x_n$ of each term 
of the sequence satisfying (\ref{form}) is given by the formula  
\beq\label{exact} 
\La_n = \left(\frac{(1-\la^{-1})\la^n -(1-\la)\la^{-n}}{\la -\la^{-1}}-1\right)\log c^{\frac{1}{d}} 
+\sum_{k=1}^{n-1} \left(\frac{\la^{n-k} -\la^{k-n}}{\la -\la^{-1}}\right)\al_k, 
\eeq 
where $\al_k$ is defined as in (\ref{lineq}) and $\la$ as in (\ref{root}).  
\end{propn} 
\begin{cor} 
To leading order, the asymptotic approximation of 
the logarithm $\La_n$ is given by (\ref{asymp}), 
where 
$$ 
C= \left(\frac{1-\la^{-1}}{\la -\la^{-1}}\right) \log c^{\frac{1}{d}} +\frac{1}{\la -\la^{-1}}\sum_{k=1}^\infty \la^{-k}\al_k, 
$$ 
and for the terms of the sequence 
$$ 
x_n \sim c^{-\frac{1}{d}}\exp (C\la^n ).  
$$  
\end{cor} 
\begin{rem} 
The form of the expression (\ref{exact}) is the discrete analogue of the solution of a linear inhomogeneous differential 
equation as obtained via the method of variation of parameters \cite{ince}; if $\al_k$ were given in advance, then 
it would provide the exact solution of (\ref{lineq}), viewed as a linear equation with initial values $\La_0=\La_1=0$.  
However, as was pointed out for a different example by Aho and Sloane \cite{ahos}, a formula such as (\ref{exact}) only represents the solution 
of the corresponding nonlinear equation (in this case the equation (\ref{form}) for $x_n$) 
in a tautologous sense, because $\al_k$ depends explicitly on the terms of the 
sequence $(x_n)$. Nevertheless, this formula does yield useful asymptotic information about the sequence.      
\end{rem}

\section{Acknowledgments} The author is grateful to Paul Hanna, Mitch Harris, Gerald McGarvey and Jeffrey Shallit 
for interesting correspondence in December 2005 related to the sequence \seqnum{A112373}; and he would like to 
thank an anonymous referee for several helpful comments, in particular  for pointing out the reference \cite{ds}.

\bigskip
\hrule
\bigskip

\noindent 2010 {\it Mathematics Subject Classification}: Primary
11J70; Secondary 11B37.

\noindent {\it Keywords}: continued fraction, nonlinear recurrence, transcendental number, Laurent property.

\bigskip
\hrule
\bigskip

\noindent (Concerned with sequences
\seqnum{A112373},
\seqnum{A114550},
\seqnum{A114551}, and
\seqnum{A114552}.)

\end{document}